\documentclass[aip,jmp,msc,preprint]{revtex4-1}
\usepackage{makeidx}
\usepackage{amsfonts}
\usepackage{amsmath}
\usepackage{amssymb}
\usepackage{eurosym}
\usepackage{graphicx}

\def\be{\begin{equation}}
\def\ee{\end{equation}}
\def\bea{\begin{eqnarray}}
\def\eea{\end{eqnarray}}

\begin{document}

\title{A new perspective on the Kosambi-Cartan-Chern theory, and its
applications}
\author{Tiberiu Harko}
\email{t.harko@ucl.ac.uk}
\affiliation{Department of Mathematics, University College London, Gower Street, London
WC1E 6BT, United Kingdom.}
\author{Praiboon Pantaragphong}
\email{kppraibo@kmitl.ac.th}
\affiliation{Mathematics Department, King Mongkut's Institute of Technology,
Ladkrabang BKK 10520, Thailand,}
\author{Sorin Sabau}
\email{sorin@tokai.ac.jp}
\affiliation{School of Science, Department of Mathematics,
Tokai University, Sapporo 005 – 8600, Japan}


\begin{abstract}
A powerful mathematical method for the investigation of the properties of
dynamical systems is represented by the Kosambi-Cartan-Chern (KCC)
theory. In this approach the time evolution of a dynamical system is described in
geometric terms, treating the solution curves of a dynamical system by geometrical methods inspired by the geodesics theory of Finsler spaces. In order to geometrize the dynamical evolution one introduces
a non-linear and a Berwald type connection, respectively, and thus the properties of any dynamical system are described in
terms of five geometrical invariants, with the second one giving the Jacobi
stability of the system. The Jacobi (in)stability is a natural
generalization of the (in)stability of the geodesic flow on a differentiable
manifold endowed with a metric (Riemannian or Finslerian) to the {\it  non-metric}
setting. Usually, the KCC theory is formulated by reducing the dynamical
evolution equations to a set of second order differential equations. In this
paper we present an alternative view on the KCC theory, in which the theory
is applied to a first order dynamical system. After introducing the general
framework of the KCC theory, we investigate in detail the properties of the
two dimensional autonomous dynamical systems. The relationship between the
linear stability and the Jacobi stability is also established. As a physical application of the formalism we consider the geometrization of Hamiltonian systems with one degree of freedom, and their stability properties.

{{\bf Keywords}: dynamical systems: curvature deviation tensor: stability: Kosambi-Cartan-Chern theory}
\end{abstract}

\maketitle

\section{Introduction}

One of the fundamental theoretical tools  extensively used in every
branch of natural sciences for modeling the evolution of natural phenomena
are the continuously time evolving dynamical systems. In scientific
applications they are extremely useful, and their usefulness is determined
by their predictive power. In turn, this predictive power is mostly
determined by the stability of their solutions. In a realistic setup some
uncertainty in the measured initial conditions in a physical system always
does exist. Therefore a mathematical model that is physically meaningful
must offer information and control on the time evolution of the deviations
of the trajectories of the dynamical system with respect to a given
reference trajectory. It is important to note that a local understanding of
the stability is equally important as the global understanding of the
late-time deviations. From a mathematical point of view the global stability
of the solutions of the dynamical systems is described by the well known
theory of the Lyapunov stability. In this approach to stability the basic
quantities are the Lyapunov exponents, measuring exponential deviations from
the given reference trajectory \cite{1,2}. However, one should mention that
it is usually very difficult to determine the Lyapunov exponents
analytically. Therefore, various numerical methods for their calculation
have been developed, and are applied in the study of the dynamical systems
\cite{3}-\cite{12}.

However, even though the methods of the Lyapunov stability analysis are well
established and well understood, it is important to adopt in the study the
stability of the dynamical system different points of view. Then, once such
a study is done, one can compare the obtained alternative results with the
corresponding Lyapunov linear stability analysis. An important alternative
approach to the study of the dynamical systems is represented by what we may
call the geometro-dynamical approach. An example of such an approach is
Kosambi-Cartan-Chern (KCC) theory, which was initiated in the pioneering
works of Kosambi \cite{Ko33}, Cartan \cite{Ca33} and Chern \cite{Ch39},
respectively. The KCC theory is inspired, and based, on the geometry of the
Finsler spaces. Its basic idea is the fundamental assumption that there is a
one to one correspondence between a second order dynamical system, and the
geodesic equations in an associated Finsler space (for a recent review of
the KCC theory see \cite{rev}). From a geometric point of view the KCC
theory is a differential geometric theory of the variational equations
describing the deviations of the whole trajectory of a dynamical system with
respect to the nearby ones \cite{An00}. In this geometrical description to
each dynamical system one associates a non-linear connection, and a Berwald
type connection. With the help of these two connections five geometrical
invariants are obtained. The most important of them is the second invariant,
also called the curvature deviation tensor. Its importance relies on the
fact that it gives the Jacobi stability of the system \cite{rev, An00,
Sa05,Sa05a}. The KCC theory has been extensively applied for the study of
different physical, biochemical or engineering systems \cite{Sa05, Sa05a,
An93, YaNa07, Ha1, Ha2, T0, T1, Ab1, Ab2, Ab3, Ab4, Ha3}.

An alternative geometrization method for dynamical systems was introduced in
\cite{Pet1} and \cite{Kau}, and further developed in \cite{Pet0}-\cite{Pet4}%
. Applications to the Henon-Heiles system and Bianchi type IX cosmological
models were also investigated. In particular, in \cite{Pet0} a theoretical
approach describing geometrically the behavior of dynamical systems, and of
their chaotic properties was considered. In this case for the underground
manifold a Finsler space was adopted. The properties of the Finsler space
allow the description of a wide class of dynamical systems, including those
with potentials depending on time and velocities. These are systems for
which the Riemannian geometry approach is generally unsuitable.

The Riemannian geometric approach to dynamical systems is based on the
well-known results that the flow associated with a time dependent
Hamiltonian
\begin{equation}
H=\frac{1}{2}\delta ^{ab}p_{a}p_{b}+V\left( x^{a}\right) ,
\end{equation}%
can be described as a geodesic flow in a curved, but conformally flat,
manifold \cite{Kau}. With the introduction of a metric of the form $%
ds^{2}=W\left( x^{a}\right) \delta _{ab}dx^{a}dx^{b},$ in which the
conformal factor is given by $W\left( x^{a}\right) =E-V\left( x^{a}\right)$ ,
where $E$ is the conserved energy associated with the time-independent
Hamiltonian $H$, it follows that in the metric $g_{ab}=W\left( x^{a}\right)
\delta _{ab},$ the geodesic equation for motion is completely equivalent to
the Hamilton equations \cite{Kau}
\begin{equation}
\frac{dx^{a}}{dt}=\frac{\partial H}{\partial p_{a}},\frac{dp_{a}}{dt}=-\frac{%
\partial H}{\partial x_{a}}.  \label{Ham}
\end{equation}%
This result implies that the confluence, or divergence of nearby
trajectories $x^{a}(s)$ and $[x+\xi ]^{a}(s)$ of the Hamiltonian dynamical
system is determined by the Jacobi equation, i.e., the equation of geodesic
deviation, which for the present case takes the following form,
\begin{equation}
\frac{D^{2}\xi ^{a}}{Ds^{2}}=R_{bcd}^{a}u^{b}u^{d}\xi ^{c}\equiv
-K_{c}^{a}\xi ^{c}.  \label{JR}
\end{equation}%
In Eq.~(\ref{JR}) $R_{bcd}^{a}$ is the Riemann tensor associated with the
metric $g_{ab}$, and $D/Ds=u^{a}\nabla _{a}$ denotes a directional
derivative along the velocity field $u^{a}=dx^{a}/ds$. Therefore linear
stability for the trajectory $x^{a}(s)$ is thus related to the Riemann
curvature $R_{bcd}^{a}$ or, more exactly, to the curvature $K_{c}^{a}$. If,
for example, $R_{bcd}^{a}$ is everywhere negative, then it follows that $%
K_{c}^{a}$ always has one or more negative eigenvalues, and therefore the
trajectory must be linearly unstable \cite{Kau}.

There are a large number of mathematical results on the geometrization of
the dynamical systems. For example, from a geometric perspective point of
view, in \cite{Punzi} the global and local stability of solutions of
continuously evolving dynamical systems was reconsidered, and the local stability was defined based on the choice of a linear connection.  Note that an important point in favor of the use of a linear
connection is that it is naturally defined for any dynamical system
$(S,X)$, and not only for those related to second-order evolution equations.

An important testing ground of the KCC theory is represented by the study of two-dimensional autonomous systems, and of their stability properties. Such a study was performed in \cite{Sa05a, Sa05} for two dimensional systems of the form
\begin{equation}\label{4n}
  \frac{du}{dt}=f(u,v),\qquad \frac{dv}{dt}=g(u,v),
\end{equation}
under the assumption that the point $(0,0)$ is a fixed point, i.e. $f(0,0)=g(0,0)=0$. By relabeling $v$ as $x$, and $g(u,v)$ as $y$, and by assuming that $g_u|_{(0,0)}\neq 0$, one can eliminate the variable $u$. Moreover, since $(u,v)=(0,0)$ is a fixed point, from the Theorem of Implicit Functions it follows that the equation $g(u,x)-y=0$ has a solution $u=u(x,y)$ in the vicinity of $(x,y)=(0,0)$. Since $\ddot x = \dot g = g_u \, f + g_v \, y$, we obtain an autonomous one-dimensional second order equation, equivalent to the system (\ref{4n}), namely
\be\label{5n}
\ddot x^1 + g^1(x,y) = 0,
\ee
where
\begin{equation}
  g^1(x,y)=-g_u(u(x,y),x) \, f(u(x,y),x) - g_v(u(x,y),x) \, y.
\end{equation}
The Jacobi stability properties of Eq.~(\ref{5n}) can be studied by using the KCC theory \cite{Sa05a,Sa05}, and the comparative study of the Jacobi and Lyapunov stability can be performed in detail.

In the present paper we will introduce an alternative view of the KCC theory by adopting the perspective of
the first order dynamical systems. After a brief presentation of the general
formalism, we will concentrate our attention to the simple, but important
case of two dimensional autonomous dynamical systems, whose properties are
studied in detail. Instead of reducing the two-dimensional autonomous system to a single second order differential equation of the form (\ref{5n}), by taking the time derivative of each equation we transform the system to an equivalent system of two second-order differential equations.
We also clarify the relation of the Jacobi stability
approach with classical Lyapunov stability theory. As a physical
application of our formalism we apply it to the study of two dimensional Hamiltonian systems, describing physical processes with one degree of freedom. The Hamilton system is transformed into an equivalent system of two second order differential equations, which can be studied geometrically similarly as geodesic equations into an associated Finsler space. We obtain the geometrical quantities describing the geometerized Hamilton system, and the conditions for its Jacobi stability are obtained.

The present paper is organized as follows. The KCC theory is presented in Section~%
\ref{kcc}. The applications of the KCC stability theory to the case of
two-dimensional systems is considered in Section~\ref{2dim}. The comparison
between the linear Lyapounov stability and the Jacobi stability is performed
in Section~\ref{comp}. The application of the KCC theory to two dimensional
Hamiltonian systems is discussed in Section~\ref{aHam}. We discuss and
conclude our results in Section~\ref{concl}.

\section{Brief review of the KCC theory and of the Jacobi stability}

\label{kcc}

In the present Section we briefly summarize the basic concepts and results
of the KCC theory, and we introduce the relevant notations (for a detailed
presentation see \cite{rev} and \cite{An00}).

\subsection{Geometrical interpretation of dynamical systems}

In the following we assume that $\mathcal{M}$ is a real, smooth $n$%
-dimensional manifold, and we denote by $T\mathcal{M}$ its tangent bundle.
On an open connected subset $\Omega $ of the Euclidian $(2n+1)$ dimensional
space $R^{n}\times R^{n}\times R^{1}$, we introduce a $2n+1$ dimensional
coordinates system $\left(x^i,y^i,t\right)$, $i=1,2,...,n$, where $\left(
x^{i}\right) =\left( x^{1},x^{2},...,x^{n}\right) $, $\left( y^{i}\right)
=\left( y^{1},y^{2},...,y^{n}\right) $ and $t$ is the usual time coordinate.
The coordinates $y^i$ are defined as
\begin{equation}
y^{i}=\left( \frac{dx^{1}}{dt},\frac{dx^{2}}{dt},...,\frac{dx^{n}}{dt}%
\right) .
\end{equation}

A basic assumption in our approach is the assumption that the time
coordinate $t$ is an absolute invariant. Therefore, the only admissible
coordinate transformations are
\begin{equation}
\tilde{t}=t,\tilde{x}^{i}=\tilde{x}^{i}\left( x^{1},x^{2},...,x^{n}\right)
,i\in \left\{1 ,2,...,n\right\} .  \label{ct}
\end{equation}

Following \cite{Punzi}, we assume that a deterministic dynamical systems can
be defined as a set of formal rules that describe the evolution of points in
a set $S$ with respect to an external time parameter $t\in T$, which can be
discrete, or continuous. More exactly, a dynamical system is a map \cite%
{Punzi}
\begin{equation}
\phi:T \times S \rightarrow S, (t,x)\mapsto \phi (t,x),
\end{equation}
which satisfies the condition $\phi (t , \cdot) \circ \phi (s , \cdot)=\phi
(t+s , \cdot)$, $\forall t ,s\in T$. For realistic dynamical systems that
can model natural phenomena additional structures need to be added to the
above definition.

In many situations of physical interest the equations of motion of a
dynamical system follow from a Lagrangian $L$ via the Euler-Lagrange
equations,
\begin{equation}
\frac{d}{dt}\frac{\partial L}{\partial y^{i}}-\frac{\partial L}{\partial
x^{i}}=F_{i},i=1,2,...,n,  \label{EL}
\end{equation}%
where $F_{i}$, $i=1,2,...,n$, is the external force. Note that the triplet $%
\left( M,L,F_{i}\right) $ is called a Finslerian mechanical system \cite%
{MiFr05,MHSS}. If the Lagrangian $L$ is regular, it follows that the
Euler-Lagrange equations defined in Eq.~(\ref{EL}) are equivalent to a
system of second-order ordinary (usually nonlinear) differential equations
\begin{equation}
\frac{d^{2}x^{i}}{dt^{2}}+2G^{i}\left( x^{j},y^{j},t\right) =0,i\in \left\{
1,2,...,n\right\} ,  \label{EM}
\end{equation}%
where each function $G^{i}\left( x^{j},y^{j},t\right) $ is $C^{\infty }$ in
a neighborhood of some initial conditions $\left( \left( x\right)
_{0},\left( y\right) _{0},t_{0}\right) $ in $\Omega $.


The fundamental idea of the KCC theory is that if an arbitrary system of
second-order differential equations of the form (\ref{EM}) is given, with no
\textit{a priori} Lagrangean function assumed, still one can study the
behavior of its trajectories by analogy with the trajectories of the
Euler-Lagrange system.

\subsection{The non-linear connection and the KCC invariants associated to a
dynamical system}

As a first step in the analysis of the geometry associated to the dynamical
system defined by Eqs.~(\ref{EM}), we introduce a nonlinear connection $N$
on $M$, with coefficients $N_{j}^{i}$, defined as \cite{MHSS}
\begin{equation}  \label{NC}
N_{j}^{i}=\frac{\partial G^{i}}{\partial y^{j}}.
\end{equation}

Geometrically the nonlinear connection $N_{j}^{i}$ can be interpreted in
terms of a dynamical covariant derivative $\nabla ^N$: for two vector fields
$v$, $w$ defined over a manifold $M$, we define the covariant derivative $%
\nabla ^N$ as \cite{Punzi}
\begin{equation}  \label{con}
\nabla _v^Nw=\left[v^j\frac{\partial }{\partial x^j}w^i+N^i_j(x,y)w^j\right]%
\frac{\partial }{\partial x^i}.
\end{equation}

For $N_i^j(x,y)=\Gamma _{il}^j(x)y^l$, from Eq.~(\ref{con}) we recover the
definition of the covariant derivative for the special case of a standard
linear connection, as defined in Riemmannian geometry.

For the non-singular coordinate transformations introduced through Eqs.~(\ref%
{ct}), the KCC-covariant differential of a vector field $\xi ^{i}(x)$ on the
open subset $\Omega \subseteq R^{n}\times R^{n}\times R^{1}$ is defined as
\cite{An93,An00,Sa05,Sa05a}
\begin{equation}
\frac{D\xi ^{i}}{dt}=\frac{d\xi ^{i}}{dt}+N_{j}^{i}\xi ^{j}.  \label{KCC}
\end{equation}

For $\xi ^{i}=y^{i}$ we obtain
\begin{equation}
\frac{Dy^{i}}{dt}=N_{j}^{i}y^{j}-2G^{i}=-\epsilon ^{i}.
\end{equation}
The contravariant vector field $\epsilon ^{i}$ defined on $\Omega $ is
called the first KCC invariant.

Now we vary the trajectories $x^{i}(t)$ of the system (\ref{EM}) into nearby
ones according to the rule
\begin{equation}
\tilde{x}^{i}\left( t\right) =x^{i}(t)+\eta \xi ^{i}(t),  \label{var}
\end{equation}
where $\left| \eta \right| $ is a small parameter, and $\xi ^{i}(t)$ are the
components of a contravariant vector field defined along the trajectory $%
x^{i}(t)$. By substituting Eqs.~(\ref{var}) into Eqs.~(\ref{EM}), and by
taking the limit $\eta \rightarrow 0$, we obtain the deviation, or Jacobi,
equations in the form \cite{An93,An00,Sa05,Sa05a}
\begin{equation}
\frac{d^{2}\xi ^{i}}{dt^{2}}+2N_{j}^{i}\frac{d\xi ^{j}}{dt}+2\frac{\partial
G^{i}}{\partial x^{j}}\xi ^{j}=0.  \label{def}
\end{equation}

Eq.~(\ref{def}) can be rewritten in a covariant form with the use of the
KCC-covariant derivative as
\begin{equation}
\frac{D^{2}\xi ^{i}}{dt^{2}}=P_{j}^{i}\xi ^{j},  \label{JE}
\end{equation}
where we have denoted
\begin{equation}  \label{Pij}
P_{j}^{i}=-2\frac{\partial G^{i}}{\partial x^{j}}-2G^{l}G_{jl}^{i}+ y^{l}%
\frac{\partial N_{j}^{i}}{\partial x^{l}}+N_{l}^{i}N_{j}^{l}+\frac{\partial
N_{j}^{i}}{\partial t},
\end{equation}
and we have introduced the Berwald connection $G_{jl}^{i}$, defined as \cite%
{rev, An00, An93,MHSS,Sa05,Sa05a}
\begin{equation}
G_{jl}^{i}\equiv \frac{\partial N_{j}^{i}}{\partial y^{l}}.
\end{equation}

The tensor $P_{j}^{i}$ is called the second KCC-invariant, or the deviation curvature
tensor, while Eq.~(\ref{JE}) is called the Jacobi equation. When the system
of equations (\ref{EM}) describes the geodesic equations, in either Riemann
or Finsler geometry Eq.~(\ref{JE}) is the Jacobi field equation.

The trace $P$ of the curvature deviation tensor can be obtained from the
relation
\begin{equation}
P=P_{i}^{i}=-2\frac{\partial G^{i}}{\partial x^{i}}-2G^{l}G_{il}^{i}+ y^{l}%
\frac{\partial N_{i}^{i}}{\partial x^{l}}+N_{l}^{i}N_{i}^{l}+\frac{\partial
N_{i}^{i}}{\partial t}.
\end{equation}

One can also introduce the third, fourth and fifth invariants of the system (%
\ref{EM}), which are defined as \cite{An00}
\begin{equation}  \label{31}
P_{jk}^{i}\equiv \frac{1}{3}\left( \frac{\partial P_{j}^{i}}{\partial y^{k}}-%
\frac{\partial P_{k}^{i}}{\partial y^{j}}\right) ,P_{jkl}^{i}\equiv \frac{%
\partial P_{jk}^{i}}{\partial y^{l}},D_{jkl}^{i}\equiv \frac{\partial
G_{jk}^{i}}{\partial y^{l}}.
\end{equation}

Geometrically, the third invariant $P_{jk}^{i}$ can be interpreted as a
torsion tensor. The fourth and fifth invariants $P_{jkl}^{i}$ and $%
D_{jkl}^{i}$ are called the Riemann-Christoffel curvature tensor, and the
Douglas tensor, respectively \cite{rev, An00}. Note that in a Berwald space
these tensors always exist. In the KCC theory they describe the geometrical
properties and interpretation of a system of second-order differential
equations.

\subsection{The definition of the Jacobi stability}

The behavior of the trajectories of the dynamical system given by Eqs.~(\ref%
{EM}) in a vicinity of a point $x^{i}\left( t_{0}\right) $ is extremely
important in many physical, chemical or biological applications. For
simplicity in the following we take $t_{0}=0$. In the following we consider
the trajectories $x^{i}=x^{i}(t)$ as curves in the Euclidean space $\left(
R^{n},\left\langle .,.\right\rangle \right) $, where $\left\langle
.,.\right\rangle $ is the canonical inner product of $R^{n}$. We assume that
the deviation vector $\xi $ obeys the initial conditions $\xi \left(
0\right) =O$ and $\dot{\xi}\left( 0\right) =W\neq O$, where $O\in R^{n}$ is
the null vector \cite{rev, An00, Sa05,Sa05a}.


Thus, for the focusing tendency of the trajectories around $t_{0}=0$ we
introduce the following description: if $\left| \left| \xi \left( t\right)
\right| \right| <t^{2}$, $t\approx 0^{+}$, then the trajectories are
bunching together. But if $\left| \left| \xi \left( t\right) \right| \right|
>t^{2}$, $t\approx 0^{+}$, the trajectories have a dispersing behavior \cite%
{rev, An00, Sa05,Sa05a}. The focusing/dispersing tendency of the
trajectories of a dynamical system can be also described in terms of the
deviation curvature tensor in the following way: The trajectories of the
system of equations (\ref{EM}) are bunching together for $t\approx 0^{+}$ if
and only if the real part of the eigenvalues of the deviation tensor $%
P_{j}^{i}\left( 0\right) $ are strictly negative. On the other hand, the
trajectories are dispersing if and only if the real part of the eigenvalues
of $P_{j}^{i}\left( 0\right) $ are strictly positive \cite{rev, An00,
Sa05,Sa05a}.

Based on the above considerations we introduce the concept of the Jacobi
stability for a dynamical system as follows \cite{rev, An00,Sa05,Sa05a}:

\textbf{Definition:} If the system of differential equations Eqs.~(\ref{EM})
satisfies the initial conditions $\left| \left| x^{i}\left( t_{0}\right) -%
\tilde{x}^{i}\left( t_{0}\right) \right| \right| =0$, $\left| \left| \dot{x}%
^{i}\left( t_{0}\right) -\tilde{x}^{i}\left( t_{0}\right) \right| \right|
\neq 0$, with respect to the norm $\left| \left| .\right| \right| $ induced
by a positive definite inner product, then the trajectories of Eqs.~(\ref{EM}%
) are Jacobi stable if and only if the real parts of the eigenvalues of the
deviation tensor $P_{j}^{i}$ are strictly negative everywhere. Otherwise,
the trajectories are Jacobi unstable.


\section{ The case of the two variable dependent dynamical system}

\label{2dim}

In the present Section we consider the case of the two-dimensional dynamical
systems already studied in \cite{Sa05,Sa05a}. However, the present
formulation is in a parametric form.

Let's consider the following two dimensional dynamical system,
\begin{equation}
\frac{dx^{1}}{dt}=f\left( x^{1},x^{2}\right) {,}  \label{ex1}
\end{equation}
\begin{equation}
\frac{dx^{2}}{dt}=g\left( x^{1},x^{2}\right) {.}  \label{ex2}
\end{equation}

We point out that we regard a solution of Eqs. (\ref{ex1}) and (\ref{ex2})
as a flow $\varphi _{t}:D\subset \mathbb{R}^{2}\rightarrow \mathbb{R} ^{2}$,
or, more generally, $\varphi _{t}:D\subset M\rightarrow M$, where $M$ is a
smooth surface in $\mathbb{R} ^{3}$. The canonical lift of $\varphi _{t}$ to
the tangent space $TM$ can be geometrically defined as
\begin{equation}
\hat{\varphi}_{t}:TM\rightarrow TM,  \label{3.2}
\end{equation}
\begin{equation}  \label{3.2a}
\hat{\varphi}_{t}(u)=\left( \varphi _{t}(u),\dot{\varphi}_{t}(u)\right) .
\end{equation}

In terms of dynamical systems, we simply take the derivative of Eqs.~(\ref%
{ex1}) and (\ref{ex2}) and obtain
\begin{equation}
\frac{d^{2}x^{1}}{dt^{2}}=f_{1}\left( x^{1},x^{2}\right) y^{1}+f_{2}\left(
x^{1},x^{2}\right) y^{2},  \label{3.3a}
\end{equation}
\begin{equation}
\frac{d^{2}x^{2}}{dt^{2}}=g_{1}\left( x^{1},x^{2}\right) y^{1}+g_{2}\left(
x^{1},x^{2}\right) y^{2},  \label{3.3b}
\end{equation}%
where we have denoted
\begin{equation}
f_{1}:=\frac{\partial f}{\partial x^{1}},f_{2}:=\frac{\partial f}{\partial
x^{2}},g_{1}:=\frac{\partial g}{\partial x^{1}},g_{2}:=\frac{\partial g}{%
\partial x^{2}},y^{1}=\frac{dx^{1}}{dt},y^{2}=\frac{dx^{2}}{dt}.
\end{equation}

In other words, on $TM$ we obtain
\begin{equation}
\frac{dy^{1}}{dt}=f_{1}\left( x^{1},x^{2}\right) y^{1}+f_{2}\left(
x^{1},x^{2}\right) y^{2},
\end{equation}
\begin{equation}
\frac{dy^{2}}{dt}=g_{1}\left( x^{1},x^{2}\right) y^{1}+g_{2}\left(
x^{1},x^{2}\right) y^{2},
\end{equation}%
where $\left( x^{1}x^{2},y^{1},y^{2}\right) $ are local coordinates on $TM$.
Hence we can see that the above system is actually a linear dynamical system
on the fiber $T_{\left( x^{1},x^{2}\right) }M$. Moreover, the system of Eqs.
(\ref{3.3a}) and (\ref{3.3b}) can be written as
\begin{equation}
\frac{d^{2}x^{1}}{dt^{2}}+\left( -f_{1}y^{1}-f_{2}y^{2}\right) =0,
\end{equation}
\begin{equation}
\frac{d^{2}x^{2}}{dt^{2}}+\left( -g_{1}y^{1}-g_{2}y^{2}\right) =0.
\end{equation}

By comparison with Eqs. (\ref{EM}) we have

\begin{equation}
\begin{pmatrix}
G^{1} \\
G^{2}%
\end{pmatrix}%
=-\frac{1}{2}%
\begin{pmatrix}
f_{1}y^{1}+f_{2}y^{2} \\
g_{1}y^{1}+g_{2}y^{2}%
\end{pmatrix}%
=-\frac{1}{2}%
\begin{pmatrix}
f_{1} & f_{2} \\
g_{1} & g_{2}%
\end{pmatrix}%
\begin{pmatrix}
y^{1} \\
y^{2}%
\end{pmatrix}%
=-\frac{1}{2}J\cdot y,
\end{equation}%
where
\begin{equation}
J=J\left( f,g\right) =%
\begin{pmatrix}
f_{1} & f_{2} \\
g_{1} & g_{2}%
\end{pmatrix}%
\end{equation}%
is the Jacobian of the dynamical system Eqs.~(\ref{ex1})-(\ref{ex2}). Using
Eq. (\ref{NC}) we obtain for the nonlinear connection
\begin{equation}
\left( N_{j}^{i}\right) _{i,j=1,2}=%
\begin{pmatrix}
N_{1}^{1} & N_{2}^{1} \\
N_{1}^{2} & N_{2}^{2}%
\end{pmatrix}%
=%
\begin{pmatrix}
\frac{\partial G^{1}}{\partial y^{1}} & \frac{\partial G^{1}}{\partial y^{2}}
\\
\frac{\partial G^{2}}{\partial y^{1}} & \frac{\partial G^{2}}{\partial y^{1}}%
\end{pmatrix}%
=-\frac{1}{2}%
\begin{pmatrix}
f_{1} & f_{2} \\
g_{1} & g_{2}%
\end{pmatrix}%
=-\frac{1}{2}J\left( f,g\right) .
\end{equation}

Therefore all the components of the Berwald connection cancel,
\begin{equation}
G_{jl}^{i}:=\frac{\partial N_{j}^{i}}{\partial y^{l}}\equiv 0.
\end{equation}

Then, for the components of the deviation curvature tensor $\left(
P_{j}^{i}\right) $, given in Eq.~(\ref{Pij}),  
we obtain
\begin{eqnarray*}
P_{1}^{1} &=&-2\left( \frac{\partial G^{1}}{\partial x^{1}}\right) +y^{1}%
\frac{\partial N_{1}^{1}}{\partial x^{1}}+y^{2}\frac{\partial N_{1}^{1}}{%
\partial x^{2}}+N_{l}^{1}N_{1}^{l}=\frac{\partial }{\partial x^{1}}\left(
f_{1}y^{1}+f_{2}y^{2}\right) -\frac{1}{2}y^{1}\frac{\partial f_{1}}{\partial
x^{1}}-\frac{1}{2}y^{2}\frac{\partial f_{1}}{\partial x^{2}}%
+N_{l}^{i}N_{j}^{l} \\
&=&\frac{1}{2}f_{11}y^{1}+\frac{1}{2}f_{12}y^{2}+N_{l}^{i}N_{j}^{l},
\end{eqnarray*}
\begin{equation*}
P_{1}^{2}=-2\left( \frac{\partial G^{2}}{\partial x^{1}}\right) +y^{1}\frac{%
\partial N_{1}^{2}}{\partial x^{1}}+y^{2}\frac{\partial N_{1}^{2}}{\partial
x^{2}}+N_{l}^{2}N_{2}^{l}=\frac{\partial }{\partial x^{1}}\left(
g_{1}y^{1}+g_{2}y^{2}\right) -\frac{1}{2}y^{1}\frac{\partial g_{1}}{\partial
x^{1}}-\frac{1}{2}y^{2}\frac{\partial g_{1}}{\partial x^{2}}%
+N_{l}^{2}N_{2}^{l},
\end{equation*}
and so on. Therefore we get
\begin{eqnarray*}
\left( P_{j}^{i}\right) &=&%
\begin{pmatrix}
P_{1}^{1} & P_{2}^{1} \\
P_{1}^{2} & P_{2}^{2}%
\end{pmatrix}%
=\frac{1}{2}%
\begin{pmatrix}
f_{11}y^{1}+f_{12}y^{2} & f_{12}y^{1}+f_{22}y^{2} \\
g_{11}y^{1}+g_{12}y^{2} & g_{12}y^{1}+g_{22}y^{2}%
\end{pmatrix}%
+\frac{1}{4}J_{l}^{i}\left( f,g\right) \times J_{j}^{l}\left( f,g\right) \\
&=&%
\begin{pmatrix}
f_{11}y^{1}+f_{12}y^{2} & g_{11}y^{1}+g_{12}y^{2} \\
f_{12}y^{1}+f_{22}y^{2} & g_{12}y^{1}+g_{22}y^{2}%
\end{pmatrix}%
^{t}+\frac{1}{4}J_{l}^{i}\left( f,g\right) \times J_{j}^{l}\left( f,g\right)
\\
&=&%
\begin{pmatrix}
\begin{pmatrix}
f_{11} & f_{12} \\
f_{21} & f_{22}%
\end{pmatrix}%
\begin{pmatrix}
y^{1} \\
y^{2}%
\end{pmatrix}%
\Bigg| &
\begin{pmatrix}
g_{11} & g_{12} \\
g_{21} & g_{22}%
\end{pmatrix}%
\begin{pmatrix}
y^{1} \\
y^{2}%
\end{pmatrix}%
\end{pmatrix}%
^{t}+\frac{1}{4}J_{l}^{i}\left( f,g\right) \times J_{j}^{l}\left( f,g\right)
.
\end{eqnarray*}

Therefore we have obtained the following

\textbf{Proposition 3.1} \textit{The curvature deviation tensor associated
to a second order dynamical system is given by
\begin{equation}
P=\frac{1}{2}%
\begin{pmatrix}
H_{f}\cdot y & H_{g}\cdot y%
\end{pmatrix}%
^{t}+\frac{1}{4}J^{2}\left( f,g\right) ,
\end{equation}
where $H_{f}=%
\begin{pmatrix}
f_{11} & f_{12} \\
f_{21} & f_{22}%
\end{pmatrix}%
$ is the Hessian of $f$, and similarly for $g$.}

\section{Lyapunov and Jacobi stability of two dimensional dynamical systems}\label{comp}

Similarly with \cite{rev}, and without losing generality, we assume $p=(0,0)$ is a
fixed point of Eqs.~(\ref{ex1}) and (\ref{ex2}). Then the Lyapunov stability
is governed by the characteristic equation
\begin{equation}
\mu ^{2}-\mathrm{tr}\;A\cdot \mu +\det A=0,  \label{p}
\end{equation}%
where tr$\;A$ and $\det A$ are the trace and determinant of the matrix
\begin{equation}
A:=J(f,g)|_{(0,0)}=%
\begin{pmatrix}
f_{1} & f_{2} \\
g_{1} & g_{2}%
\end{pmatrix}%
|_{(0,0)}.
\end{equation}

We also denote by $\Delta =\left( \mathrm{tr}\;A\right) ^{2}-4\det A$ the
discriminant of Eq. (\ref{p}). By descending again on $\mathbb{R} ^{2}$ (or $%
M$), and evaluate at the fixed point $\left( 0,0\right) $, we obtain $\left(
y^{1},y^{2}\right) |_{(0,0)}=(0,0)$, and therefore
\begin{equation}
P|_{(0,0)}=\left( \frac{1}{2}A\right) ^{2}.
\end{equation}

To be more precise, we have
\begin{equation}
P|_{(0,0)}=\frac{1}{4}%
\begin{pmatrix}
f_{1} & f_{2} \\
g_{1} & g_{2}%
\end{pmatrix}%
\begin{pmatrix}
f_{1} & f_{2} \\
g_{1} & g_{2}%
\end{pmatrix}%
=\frac{1}{4}%
\begin{pmatrix}
f_{1}^{2}+f_{2}g_{1} & f_{1}f_{2}+f_{2}g_{2} \\
f_{1}g_{1}+g_{1}g_{2} & f_{2}g_{1}+g_{2}^{2}%
\end{pmatrix}%
.
\end{equation}

Hence we have
\begin{equation}
\mathrm{tr}\;P=\frac{1}{4}\left( f_{1}^{2}+2f_{2}g_{1}+g_{2}^{2}\right) =%
\frac{1}{4}\left(
f_{1}^{2}+2f_{1}g_{2}+g_{2}^{2}-2f_{1}g_{2}+2f_{2}g_{1}\right) =\frac{1}{4}%
\left[ \left( \mathrm{tr}\;A\right) ^{2}-2\det A\right] ,
\end{equation}
\begin{equation}
\det P=\frac{1}{4}\left( f_{1}g_{2}-f_{2}g_{1}\right) ^{2}=\left(\frac{1}{4}
\det A\right) ^{2}.
\end{equation}

The eigenvalues of the matrix $\left( P_{j}^{i}\right) $ are given by the
characteristic equation
\begin{equation}
\lambda ^{2}-\mathrm{tr}\;P\cdot \lambda +\det P=0,
\end{equation}
and its discriminant is
\begin{equation}
\tilde{\Delta}=\left( {\rm tr}\;P\right) ^{2}-4\det P=\frac{1}{16}\left[ \left(
{\rm tr}\;A\right) ^{2}-2\det A\right] ^{2}-\frac{1}{4}\left( \det A\right) ^{2}=%
\frac{1}{16}\left( {\rm tr}\;A\right) ^{2}\left[ \left( {\rm tr }\;A\right) ^{2}-4\det A%
\right] .
\end{equation}

Thus we obtain the following

\textbf{Computational Lemma 4.1.} \textit{The trace, determinant and
discriminant of the characteristic equation of the deviation curvature
matrix $P$ are}:
\begin{equation*}
\mathrm{tr}\;P=\frac{1}{4}\left[ \left( {\rm tr}\;A\right) ^{2}-2\det A\right] ,
\end{equation*}
\begin{equation*}
\det P=\frac{1}{16}\left( \det A\right) ^{2},
\end{equation*}
\begin{equation*}
\tilde{\Delta }=\frac{1}{16}\left({\rm tr}\;A\right) ^{2}\Delta ,
\end{equation*}
\textit{where $\Delta =\left({\rm tr}\;A\right)^2-4\det A$ is the discriminant of Eq.~(\ref{p}).}

We recall now some general results of linear algebra.

{\bf Lemma 4.2.} Let $A$ be a $\left( 2,2\right) $ matrix, and denote by $\lambda
_{1}$, $\lambda _{2}$ its eigenvalues. Then

(i) the eigenvalues of $k\cdot A$ are $k\cdot \lambda _{1}$, $k\cdot \lambda
_{2}$, for $k\neq 0$ scalar.

(ii) the eigenvalues of $A^{k}:=\underbrace{A\cdot A...\cdot A}$ are $\left(
\lambda _{1}\right) ^{k}$, $\left( \lambda _{2}\right) ^{k}$.

From here it follows:

\textbf{Lemma 4.3.} \textit{If $\lambda _{1}$, $\lambda _{2}$ are
eigenvalues of $A$, then $\mu _{1}=\left( \frac{1}{2}\lambda _{1}\right)
^{2} $, $\mu _{2}=\left( \frac{1}{2}\lambda _{2}\right) ^{2}$ are
eigenvalues of $P$.}

Remark: The formulas of Lemma 4.3 imply
\begin{equation*}
S:=\mu _{1}+\mu _{2}=\frac{1}{4}\left( \lambda _{1}^{2}+\lambda
_{2}^{2}\right) =\frac{1}{4}\left[ \left( \lambda _{1}+\lambda _{2}\right)
^{2}-2\cdot \lambda _{1}\lambda _{2}\right] ,
\end{equation*}

\begin{equation*}
\rho :=\mu _{1}\cdot \mu _{2}=\left( \frac{1}{2}\lambda _{1}^{2}\right)
\cdot \left( \frac{1}{2}\lambda _{2}\right) ^{2}=\frac{1}{16}\left( \lambda
_{1}\lambda _{2}\right) ^{2}.
\end{equation*}

The above formulas are consistent with the computational Lemma 4.1.

\subsection{Comparison of Jacobi and Lyapunov stability}

\textbf{Definition} Let $p$ be a fixed point of the two-dimensional system $%
\dot{\mathbf{x}}=f(\mathbf{x})$, and denote by $\lambda_1$, $\lambda_2$ the
two eigenvalues of $A := (Df)_{|p}$. The following classification of the
fixed point $p$ is standard.

\begin{itemize}
\item[(I)] $\Delta >0$, \textbf{$\lambda_1$, $\lambda_2$ are real and
distinct.}

\begin{itemize}
\item[I.1.] \textbf{$\lambda_1\cdot \lambda_2>0$ (the eigenvalues have the
same sign)}: $p$ is called a \textit{node} or type I singularity; that is,
every orbit tends to the origin in a definite direction as $t \to \infty$.

\begin{itemize}
\item[I.1.1.] \textbf{$\lambda_1, \lambda_2>0$ }: $p$ is an \textit{unstable
node}.%
\index{unstable node}

\item[I.1.2.] \textbf{$\lambda_1, \lambda_2<0$ }: $p$ is a \textit{stable
node}.%
\index{stable node}
\end{itemize}

\item[I.2.] $\Delta <0$, \textbf{$\lambda_1\cdot \lambda_2<0$ (the
eigenvalues have different signs)}: $p$ is an \textit{unstable fixed point},
or a \textit{saddle} point singularity.
\end{itemize}

\item[(II)] \textbf{$\lambda_1$, $\lambda_2$ are complex, i.e. $%
\lambda_{1,2}=\alpha\pm i \beta$, $\beta \neq 0$.}

\begin{itemize}
\item[II.1.] \textbf{$\alpha\neq 0$}: $p$ is a \textit{spiral}, or a \textit{%
focus}, that is, the solutions approach the origin as $t\to \infty$, but not
from a definite direction.

\begin{itemize}
\item[II.1.1.] \textbf{$\alpha<0$}: $p$ is a \textit{stable focus}.

\item[II.1.2.] \textbf{$\alpha>0$}: $p$ is an \textit{unstable focus}.
\end{itemize}

\item[II.2.] \textbf{$\alpha = 0$}: $p$ is a \textit{center},
\index{center} that means it is not stable in the usual sense, and we have
to look at higher order derivatives.
\end{itemize}

\item[(III)] $\Delta =0$, \textbf{$\lambda_1$, $\lambda_2$ are equal, i.e. $%
\lambda_1=\lambda_2=\lambda$.}

\begin{itemize}
\item[III.1.] If there are two linearly independent eigenvectors, we have a
\textit{star singularity}, or a \textit{stable singular node} (these are
simple straight lines through the origin).

\item[III.2.] If there is only one linearly independent eigenvector, we have
an \textit{improper node}, or \textit{unstable degenerate node}.
\end{itemize}
\end{itemize}

By combining the above results with the Computational Lemma 4.1, it follows
that

I. $%
\tilde{\Delta}>0$, $\mu _1, \mu _2 \in \mathbb{R}$.

I.1. $S>0$, $p$ is Jacobi unstable $\Leftrightarrow$
\begin{equation}  \label{3.4}
\left(\mathrm{tr}\;A\right)^2-4\det A>0, \left(\mathrm{tr}\;A\right)^2-2\det
A>0.
\end{equation}

- for $\det A>0$, we must have $\left(\mathrm{tr}\;A\right)^2>2\det A>4\det
A $, and hence both cases $\left(\mathrm{tr}\;A>0,\det A>0\right)$ and $%
\left(\mathrm{tr}\;A<0,\det A<0\right)$ imply Jacobi stability.

-for $\det A<0$, we obtain again Jacobi stability because now relations (\ref%
{3.4}) are identically satisfied.

I.2. $S>0$: this case is not possible algebraically.

II. $\tilde{\Delta}<0$: $\mu _{1,2}=\alpha\pm i\beta \in
\mathbb{C}$ (Remark that $\tilde{\Delta}<0 \Rightarrow \Delta <0$, i.e. $%
\lambda _1, \lambda _2$ are complex, and hence $\det A>0$.

II.1. $S>0$: $p$ is Jacobi unstable $\Leftrightarrow$
\begin{equation}  \label{3.41}
\left(\mathrm{tr}\;A\right)^2-4\det A<0, \left(\mathrm{tr}\;A\right)^2-2\det
A>0,
\end{equation}
that is, $2\det A<\left(\mathrm{tr }A\right)^2<4\det A$. This case is
possible for $\left(\mathrm{tr}\;A>0,\det A>0\right)$, and $\left(\mathrm{tr}%
\;A<0,\det A>0\right)$.

For $\det A>0$ this case is not possible algebraically.

II.2. $S<0$: $p$ is Jacobi stable $\Leftrightarrow$
\begin{equation}  \label{3.42}
\left(\mathrm{tr}\;A\right)^2-4\det A<0, \left(\mathrm{tr}\;A\right)^2-2\det
A<0,
\end{equation}
that is $\left(\mathrm{tr}\;A\right)^2<2\det A<4\det A$. This is possible
for both $\left(\mathrm{tr}\;A>0,\det A>0\right)$, and $\left(\mathrm{tr}%
\;A<0,\det A>0\right)$.

III. $\tilde{\Delta}=0$: $\mu _1=\mu _2\in \mathbb{R}$, equivalent to $%
\Delta =0$, or $\mathrm{tr}\;A=0$.

- if $\Delta =0$, then $\lambda _1=\lambda _2\in \mathbb{R}$ and $p$ is
singular node or degenerate node

-if $\mathrm{tr}\;A=0$, $\Delta \neq 0$, it follows

(i) $\Delta >0$, $\det A<0$, i.e. saddle point, or

(ii) $\Delta <0$, $\det A<0$, i.e., $p$ is a center.

\subsection{Important remarks}

1. Assume $p$ is Jacobi stable, that is, by definition, we must have one of
the following situations:

(i) $\tilde{\Delta }>0$, $S>0$, or

(ii) $\tilde{\Delta }<0$, $S>0$.

As we have seen already (i) is algebraically impossible, hence, if $p$ is
Jacobi stable it must follow $\tilde{\Delta }<0 \Leftrightarrow \Delta <0$.
Hence we have proved:

If $p$ is Jacobi stable, then $\Delta <0$.

2. Conversely, assume $\Delta <0$. This is equivalent to $\tilde {\Delta }<0$
due to the Computational Lemma 3.3. Since $\Delta <0$ implies $\det A>0$,
and $\lambda _{1,2}=\alpha \pm i\beta $, we have
\begin{equation}
\mu _1=\left(\frac{1}{2}\lambda _1\right)^2=\frac{1}{4}\left(\alpha +i\beta
\right)^2=\frac{1}{4}\left[\left(\alpha ^2-\beta ^2\right)+2i\alpha \beta%
\right],
\end{equation}
\begin{equation}
\mu _2=\left(\frac{1}{2}\lambda _2\right)^2=\frac{1}{4}\left(\alpha -i\beta
\right)^2=\frac{1}{4}\left[\left(\alpha ^2-\beta ^2\right)-2i\alpha \beta %
\right].
\end{equation}

From the above equations we obtain
\begin{equation}
S=\frac{1}{2}\left(\alpha ^2-\beta ^2\right).
\end{equation}

Therefore it follows

(i) $\alpha ^2-\beta ^2>0$ $\Rightarrow$ $p$ is Jacobi unstable,

(ii) $\alpha ^2-\beta ^2<0$ $\Rightarrow$ $p$ is Jacobi stable.

Condition (ii) above does not appear in \cite{rev} due to the reduction of
the two dimensional dynamical system to a one-dimensional SODE. We obtain

\textbf{Theorem 4.1.} 1. If $p$ is a Jacobi stable fixed point, then $\Delta
<0$.

2. If $\Delta <0$ for the fixed point $p$, then

(i) if $\alpha ^2-\beta ^2>0$, then $p$ is Jacobi unstable;

(ii) if $\alpha ^2-\beta ^2<0$, then $p$ is Jacobi stable, where $\lambda
_{1,2}=\alpha \pm i\beta $ are the eigenvalues of $A$.

\section{Applications to two dimensional Hamiltonian systems}

\label{aHam}

As a physical application of the formalism developed in this paper, in the
present Section we consider the geometrical description, and the stability
properties of a physical system described by a Hamiltonian function $%
H=H(x,p)=H\left(x^1,x^2\right)$, where $H(x,p)\in C^n$, $n\geq 2$.
From a physical point of view $x$ represents the particle's coordinate,
while $p$ is its momentum. The motion of the system is described by a
Hamiltonian system of equations in the plane.

\textbf{Definition.} A system of differential equations on $\mathbf{R}^{2}$
is called a conservative Hamiltonian system with one degree of freedom if it
can be expressed in the form
\begin{equation}
\frac{dx^{1}}{dt}=\frac{\partial H\left( x^{1},x^{2}\right) }{\partial x^{2}}%
=H_{2}\left( x^{1},x^{2}\right) ,\frac{dx^{2}}{dt}=-\frac{\partial H\left(
x^{1},x^{2}\right) }{\partial x^{1}}=-H_{1}\left( x^{1},x^{2}\right) .
\label{54}
\end{equation}

By taking the total derivative of the Hamiltonian function, and with the use
of Eqs.~(\ref{54}) we obtain
\begin{equation}  \label{55}
\frac{dH}{dt}=\frac{\partial H\left(x^1,x^2\right)}{\partial x^1}\frac{dx^1}{%
dt}+\frac{\partial H\left(x^1,x^2\right)}{\partial x^2}\frac{dx^2}{dt}\equiv
0.
\end{equation}
Eqs.~(\ref{55}) shows that $H\left(x^1,x^2\right)$ is constant along the
solution curves of Eqs.~(\ref{54}). Hence the Hamiltonian function is a
first integral and a constant of motion. Moreover, all the trajectories of
the dynamical system lie on the contours defined by $H\left(x^1,x^2\right) =
C$, where $C$ is a constant. From a physical point of view $%
H\left(x^1,x^2\right)$ represents the total energy, which is a conserved
quantity.

By taking the derivatives of Eqs.~(\ref{54}) with respect to the time
parameter $t$ we obtain first
\begin{equation}
\frac{d^{2}x^{1}}{dt^{2}}=H_{21}\left( x^{1},x^{2}\right) y^{1}+H_{22}\left(
x^{1},x^{2}\right) y^{2},
\end{equation}%
\begin{equation}
\frac{d^{2}x^{2}}{dt^{2}}=-H_{11}\left( x^{1},x^{2}\right)
y^{1}-H_{12}\left( x^{1},x^{2}\right) y^{2},
\end{equation}%
which can be written in the equivalent form
\begin{equation}
\frac{d^{2}x^{i}}{dt^{2}}+2G^{i}\left( x^{1},x^{2},y^{1},y^{2}\right)
=0,i=1,2,  \label{56}
\end{equation}%
where

\begin{equation}
\begin{pmatrix}
G^{1} \\
G^{2}%
\end{pmatrix}%
=-\frac{1}{2}%
\begin{pmatrix}
H_{21} & H_{22} \\
-H_{11} & -H_{12}%
\end{pmatrix}%
\begin{pmatrix}
y^{1} \\
y^{2}%
\end{pmatrix}%
=-\frac{1}{2}J_{H}\cdot y,
\end{equation}%
where

\begin{equation}
J_{H}=J_H\left( -H_{1},H_{2}\right) =%
\begin{pmatrix}
H_{21} & H_{22} \\
-H_{11} & -H_{12}%
\end{pmatrix}%
,
\end{equation}%
is the Jacobian of the Hamiltonian system given by Eqs.~(\ref{54}). Eqs.~(%
\ref{56}) give the geometrical interpretation of a bidimensional (one degree
of freedom) Hamiltonian system, showing that they can be studied by similar methods as the
geodesics in a Finsler space.

Using Eq. (\ref{NC}) we obtain for the nonlinear connection associated to a
Hamiltonian system the expressions

\begin{equation}
\left( N_{j}^{i}\right) _{i,j=1,2}=%
\begin{pmatrix}
N_{1}^{1} & N_{2}^{1} \\
N_{1}^{2} & N_{2}^{2}%
\end{pmatrix}%
=-\frac{1}{2}%
\begin{pmatrix}
H_{21} & H_{22} \\
-H_{11} & -H_{12}%
\end{pmatrix}%
=-\frac{1}{2}J_{H}\left( -H_{1},H_{2}\right) .
\end{equation}

Therefore, for a Hamiltonian system all the components of the Berwald
connection vanish,

\begin{equation}
G_{jl}^{i}:=\frac{\partial N_{j}^{i}}{\partial y^{l}}\equiv 0.
\end{equation}

The components of the deviation curvature tensor of a Hamiltonian system can
be obtained as

\begin{eqnarray*}
\left( P_{j}^{i}\right) _{H} &=&%
\begin{pmatrix}
P_{1}^{1} & P_{2}^{1} \\
P_{1}^{2} & P_{2}^{2}%
\end{pmatrix}%
_H=
\begin{pmatrix}
H_{211}y^{1}+H_{212}y^{2} & H_{212}y^{1}+H_{222}y^{2} \\
-H_{111}y^{1}-H_{112}y^{2} & -H_{112}y^{1}-H_{122}y^{2}%
\end{pmatrix}%
^{t}+ \\
&&\frac{1}{4}\left( J_{H}\right) _{l}^{i}\left( -H_{1},H_{2}\right) \times
\left( J_{H}\right) _{j}^{l}\left( -H_{1},H_{2}\right) \\
&=&%
\begin{pmatrix}
\begin{pmatrix}
H_{211} & H_{212} \\
-H_{111} & -H_{112}%
\end{pmatrix}%
\begin{pmatrix}
y^{1} \\
y^{2}%
\end{pmatrix}%
\Bigg| &
\begin{pmatrix}
H_{212} & H_{222} \\
-H_{112} & -H_{122}%
\end{pmatrix}%
\begin{pmatrix}
y^{1} \\
y^{2}%
\end{pmatrix}%
\end{pmatrix}%
^{t}+ \\
&&\frac{1}{4}\left( J_{H}\right) _{l}^{i}\left( -H_{1},H_{2}\right) \times
\left( J_{H}\right) _{j}^{l}\left( -H_{1},H_{2}\right) .
\end{eqnarray*}

Therefore we have the following

\textbf{Proposition 5.1} \textit{The curvature deviation tensor associated
to a Hamiltonian dynamical system is given by
\begin{equation}
P_{H}=\frac{1}{2}%
\begin{pmatrix}
\mathit{H}_{-H_{1}}\cdot y & \mathit{H}_{H_{2}}\cdot y%
\end{pmatrix}%
^{t}+\frac{1}{4}J_{H}^{2}\left( -H_{1},H_{2}\right) ,
\end{equation}%
where $\mathit{H}_{-H_{1}}=%
\begin{pmatrix}
H_{211} & H_{212} \\
-H_{111} & -H_{112}%
\end{pmatrix}%
$ is the Hessian of $H_{1}$, and similarly for $H_{2}$.}

Hence in the following we can introduce the concept of the Jacobi stability
of a Hamiltonian system by means of the following

\textbf{Definition. }If the Hamiltonian system of Eqs.~(\ref{54}) satisfies
the initial conditions $\left\vert \left\vert x^{i}\left( t_{0}\right) -%
\tilde{x}^{i}\left( t_{0}\right) \right\vert \right\vert =0$, $\left\vert
\left\vert \dot{x}^{i}\left( t_{0}\right) -\tilde{x}^{i}\left( t_{0}\right)
\right\vert \right\vert \neq 0$, with respect to the norm $\left\vert
\left\vert .\right\vert \right\vert $ induced by a positive definite inner
product, then the trajectories of Hamiltonian dynamical system are Jacobi
stable if and only if the real parts of the eigenvalues of the curvature
deviation tensor $P_{H}$ are strictly negative everywhere. Otherwise, the
trajectories are Jacobi unstable.

To illustrate the implications of the geometric approach introduced here we
consider the simple case of the one dimensional conservative motion of a
point particle with mass $m>0$ under the influence of an external potential $%
V(x)$, described by the Hamiltonian
\begin{equation}  \label{Ham}
H=\frac{p^2}{2m}+V(x)=\frac{1}{2m}\left(x^2\right)^2+V\left(x^1\right).
\end{equation}
The equation of motion of the particle are
\begin{equation}  \label{65}
\frac{dx^1}{dt}=\frac{1}{m}x^2,\frac{dx^2}{dt}=-V^{\prime }\left(x^1\right),
\end{equation}
where in the following we denote by a prime the derivative with respect to
the coordinate $x^1$. By taking the derivative with respect to time of Eqs.~(%
\ref{65}) we obtain first
\begin{equation}  \label{66}
\frac{d^2x^1}{dt^2}=\frac{1}{m}\frac{dx^2}{dt}=\frac{1}{m}y^2, \frac{d^2x^2}{%
dt^2}=-V^{\prime \prime }\left(x^1\right)\frac{dx^1}{dt}=-V^{\prime \prime
}\left(x^1\right)y^1.
\end{equation}
Eqs.~(\ref{66}) can be written as
\begin{equation}  \label{67}
\frac{d^2x^i}{dt^2}+2G^i\left(x^1,x^2,y^1,y^2\right)=0,i=1,2,
\end{equation}
where
\begin{equation}
G^1\left(x^1,x^2,y^1,y^2\right)=-\frac{1}{2m}y^2,G^2\left(x^1,x^2,y^1,y^2%
\right)=\frac{1}{2}V^{\prime \prime }\left(x^1\right)y^1.
\end{equation}
The components of the non-linear connection $N^i_j=\partial G^i/\partial y^j$
can be obtained as
\begin{equation}
N_1^1=\frac{\partial G^1}{\partial y^1}=0, N^1_2=\frac{\partial G^1}{%
\partial y^2}=-\frac{1}{2m},N^2_1=\frac{\partial G^2}{\partial y^1}=\frac{1}{%
2}V^{\prime \prime }\left(x^1\right), N_2^2=\frac{\partial G^2}{\partial y^2}%
=0.
\end{equation}

With the use of Eq.~(\ref{Pij}) we obtain for the components of the
deviation curvature tensor $P^i_j$ the expressions
\begin{equation}
P_1^1=-\frac{1}{4m}V^{\prime \prime }\left(x^1\right),P^1_2=0,P^2_1=-\frac{1%
}{2}V^{\prime \prime }\left(x^1\right)y^1,P_2^2=-\frac{1}{4m}V^{\prime
\prime }\left(x^1\right).
\end{equation}
The eigenvalues $\lambda _{1,2}$ of the curvature deviation tensor are given
by
\begin{equation}
\lambda _{1,2}=\frac{1}{2}\left[P_1^1+P_2^2\pm\sqrt{4P^1_2P^2_1+%
\left(P_1^1-P_2^2\right)^2}\right],
\end{equation}
and can be obtained explicitly as
\begin{equation}
\lambda_{1,2}=-\frac{1}{4m}V^{\prime \prime }\left(x^1\right).
\end{equation}

Therefore we have obtained the following

\textbf{Theorem 5.1} The trajectories of a one-dimensional Hamiltonian
dynamical system with point particle like Hamiltonian given by Eq.~(\ref{Ham}%
) are Jacobi stable if and only if the potential $V(x)$ of the external
forces satisfies for all $x$ the condition $V^{\prime \prime }(x)>0$.

From a physical point of view the condition $V^{\prime \prime }(x)>0$
implies that the potential $V(x)$ is in a minimum. More exactly, if $%
V^{\prime }\left(x_0\right)=0$ is an equilibrium state of a physical system
at point $x_0$, the condition for the Jacobi stability $V^{\prime \prime
}\left(x_0\right)>0$ requires that the potential $V$ has a minimum at $x=x_0$%
.

For a point particle like two dimensional Hamiltonian system with one degree
of freedom the equations of the geodesic deviation, satisfied by the
deviation vector $\xi ^i$, $i=1,2$, given by Eq.~(\ref{def}), can be
obtained as
\begin{equation}
\frac{d^2\xi ^1}{dt^2}-\frac{1}{m}\frac{d\xi ^2}{dt}=0,
\end{equation}
\begin{equation}
\frac{d^2\xi ^2}{dt^2}+V^{\prime \prime }\left(x^1\right)\frac{d\xi ^1}{dt}%
+V^{\prime \prime \prime }\left(x^1\right)y^1\xi ^1=0.
\end{equation}

Let's assume now that $\left( x^{1}=x_{0},x^{2}=0,y^{1}=0,y^{2}=0\right) $
are the critical points of the Hamiltonian system of Eqs.~(\ref{65}). Then
the system of the geodesic deviation equations take the form
\begin{equation}
\frac{d^{2}\xi ^{1}}{dt^{2}}-\frac{1}{m}\frac{d\xi ^{2}}{dt},
\end{equation}%
\begin{equation}
\frac{d^{2}\xi ^{2}}{dt^{2}}+V^{\prime \prime }\left( x_{0}\right) \frac{%
d\xi ^{1}}{dt}=0,
\end{equation}%
and must be integrated with the initial conditions $\xi ^{1}(0)=0$, $\xi
^{2}(0)=0$, $\dot{\xi}^{1}=\xi _{10}$, and $\dot{\xi}^{2}(0)=\xi _{20}$,
respectively. The geodesic deviation equations for the point mass
Hamiltonian have the general solution
\begin{equation}
\xi ^{1}(t)=\frac{1}{V^{\prime \prime }\left( x_{0}\right) }\left[ \sqrt{m}%
\xi _{10}\sqrt{V^{\prime \prime }\left( x_{0}\right) }\sin \left( \frac{%
\sqrt{V^{\prime \prime }\left( x_{0}\right) }}{\sqrt{m}}t\right) -\xi
_{20}\cos \left( \frac{\sqrt{V^{\prime \prime }\left( x_{0}\right) }}{\sqrt{m%
}}t\right) +\xi _{20}\right] ,
\end{equation}%
\begin{equation}
\xi ^{2}(t)=m\xi _{10}\left[ \cos \left( \frac{\sqrt{V^{\prime \prime
}\left( x_{0}\right) }}{\sqrt{m}}t\right) -1\right] +\frac{\sqrt{m}\xi _{20}%
}{\sqrt{V^{\prime \prime }\left( x_{0}\right) }}\sin \left( \frac{\sqrt{%
V^{\prime \prime }\left( x_{0}\right) }}{\sqrt{m}}t\right) .
\end{equation}

\section{Discussions and final remarks}\label{concl}

In the present paper we have considered an alternative approach to the standard KCC theory for first order autonomous dynamical systems, based on a different transformation of the system to second order differential equations. The approach was presented in detail for two dimensional dynamical systems, for which the  two basic stability analysis methods -- the (Lyapunov) linear stability analysis and the Jacobi stability analysis -- were discussed in detail. From the point of view of the KCC theory the present approach allows an extension of the geometric framework for first order systems, such increasing the parameter space, and the predictive power, of the method. We have also found that there is a good correlation between the linear stability of the critical points, and the Jacobi stability of the same points, describing the robustness of the corresponding trajectory to a small perturbation \cite{Sa05}. On the other hand, the Jacobi stability is a very convenient way of describing the resistance of limit cycles to small perturbation of trajectories.

As an application of our approach we have considered the study of a bi-dimensional Hamiltonian system, describing the one dimensional (one degree of freedom) motion of a physical system. The KCC theory can provide an alternative, and very powerful, method for the geometrization of classical mechanical systems, whose properties can be described by a Hamiltonian function. The transformation of the corresponding Hamilton equations to second order differential equations allows naturally their study similarly as geodesics in an associated Finsler space, and gives the possibility of a full geometric description of the properties of the dynamical system in a {\it non-metric} setting. This represents one of the basic differences, and advantages, of the KCC approach as compared to the alternative Jacobi and Eisenhart methods for geometrization \cite{PR}, which essentially require as a starting point a metric. It is important to emphasize that the advantages of the geometric approach to the description of the dynamical systems are not only conceptual, but also the method has a predictive value. By starting from the deviation curvature tensor we can obtain some effective stability conditions for physical systems. Moreover, in the present approach, the geodesic deviation equation can be formulated and solved rather easily (either analytically or numerically), and thus the behavior of the full perturbations of the trajectories near critical points can be studied in detail.

To summarize our results, in the present paper we have introduced and studied in detail some geometrical theoretical tools necessary for an in depth analysis and description of the stability properties of dynamical systems that may play a fundamental role in our understanding of the evolution of natural phenomena.

\section*{Acknowledgments}


\begin{thebibliography}{99}
\bibitem{1} A. M. Mancho, D. Small, and S. Wiggins, Physics Reports \textbf{%
437}, 55 (2006).

\bibitem{2} G. Boffetta, M. Cencini, M. Falcioni, and A. Vulpiani, Physics
Reports \textbf{356}, 367 (2002).

\bibitem{3} A. E. Motter, M. Gruiz, G. K\'{a}rolyi, and T. T\'{e}l, Phys.
Rev. Lett. \textbf{111}, 194101 (2013).

\bibitem{4} E. G. Altmann, J. S. E. Portela, and T. T\'{e}l, Phys. Rev.
Lett. \textbf{111}, 144101 (2013).

\bibitem{5} C. Skokos, I. Gkolias, and S. Flach, Phys. Rev. Lett. \textbf{111%
}, 064101 (2013).

\bibitem{6} D. Paz\'{o}, J. M. L\'{o}pez, and A. Politi, Phys. Rev. \textbf{%
E 87}, 062909 (2013).

\bibitem{7} B. Schalge, R. Blender, J. Wouters, K. Fraedrich, and F.
Lunkeit, Phys. Rev. \textbf{E 87}, 052113 (2013).

\bibitem{8} S. Zeeb, T. Dahms, V. Flunkert, E. Sch\"{o}ll, I. Kanter, and W.
Kinzel, Phys. Rev. \textbf{E 87}, 042910 (2013).

\bibitem{9} C. J. Yang, W. D. Zhu, and G. X. Ren, Communications in
Nonlinear Science and Numerical Simulation \textbf{18}, 3271 (2013).

\bibitem{10} H.-Liu Yang, G. Radons, and H. Kantz, Phys. Rev. Lett. \textbf{%
109}, 244101 (2012).

\bibitem{11} A. Politi, F. Ginelli, S. Yanchuk, and Y. Maistrenko, Physica
D: Nonlinear Phenomena \textbf{224}, 90 (2006).

\bibitem{12} L. Donetti, P. I. Hurtado, and M. A. Munoz, Phys. Rev. Lett.
\textbf{95}, 188701 (2005).

\bibitem{Ko33} D. D. Kosambi, Math. Z. \textbf{37}, 608 (1933).

\bibitem{Ca33} E. Cartan, Math. Z. \textbf{37}, 619 (1933).

\bibitem{Ch39} S. S. Chern, Bulletin des Sciences Mathematiques \textbf{63},
206 (1939).

\bibitem{rev} C. G. Boehmer, T. Harko, and S. V. Sabau, Adv. Theor. Math.
Phys. 16 (2012) 1145-1196.

\bibitem{An00} P. L. Antonelli (Editor), Handbook of Finsler geometry, vol.
1, Kluwer Academic, Dordrecht, (2003).

\bibitem{Sa05a} S. V. Sabau, Nonlinear Analysis: Real World Applications
\textbf{6}, 563 (2005).

\bibitem{Sa05} S. V. Sabau, Nonlinear Analysis \textbf{63}, 143 (2005).

\bibitem{An93} P. L. Antonelli, Tensor, N. S. \textbf{52}, 27 (1993).

\bibitem{YaNa07} T. Yajima and H. Nagahama, J. Phys. A: Math. Theor. \textbf{%
40}, 2755 (2007).

\bibitem{Ha1} T. Harko and V. S. Sabau, Phys. Rev. \textbf{D 77}, 104009
(2008).

\bibitem{Ha2} C. G. Boehmer and T. Harko, Journal of Nonlinear Mathematical
Physics \textbf{17}, 503 (2010).

\bibitem{T0} T. Yajima and H. Nagahama, Acta Mathematica Academiae
Paedagogicae Ny\'{\i}regyh\'{a}ziensis \textbf{24}, 179 (2008).

\bibitem{T1} T. Yajima and H. Nagahama, Physics Letters \textbf{A 374}, 1315
(2010).

\bibitem{Ab1} H. Abolghasem, Journal of Dynamical Systems and Geometric
Theories \textbf{10}, 13 (2012).

\bibitem{Ab2} H. Abolghasem, Journal of Dynamical Systems and Geometric
Theories \textbf{10}, 197 (2012).

\bibitem{Ab3} H. Abolghasem, International Journal of Differential Equations
and Applications \textbf{12}, 131 (2013).

\bibitem{Ab4} H. Abolghasem, International Journal of Pure and Applied
Mathematics \textbf{87}, 181 (2013).

\bibitem{Ha3} T. Harko, C. Y. Ho, C. S. Leung, and S. Yip, Int. J. of
Geometric Methods in Modern Physics \textbf{12}, 1550081 (2015).

\bibitem{Pet1} M. Pettini, Phys. Rev. \textbf{E 47}, 828 (1993).

\bibitem{Kau} H. E. Kandrup, Phys. Rev. \textbf{E 56}, 2722 (1997).

\bibitem{Pet0} M. Di Bari, D. Boccaletti, P. Cipriani, and G. Pucacco, Phys.
Rev. \textbf{E 55}, 6448 (1997).

\bibitem{Pet2} P. Cipriani and M. Di Bari, Phys. Rev. Lett. \textbf{81},
5532 (1998).

\bibitem{Pet3} M. Di Bari and P. Cipriani, Planet. Space. Science \textbf{46}%
, 1543 (1998).

\bibitem{PR} L. Casetti, M. Pettini, and E. G. D. Cohen, Physics Reports
\textbf{337}, 237 (2000).

\bibitem{Pet4} G. Ciraolo and M. Pettini, Celestial Mechanics and Dynamical
Astronomy \textbf{83}, 171 (2002).

\bibitem{Punzi} R. Punzi and M. N. R. Wohlfarth, Phys. Rev. \textbf{E 79},
046606 (2009).

\bibitem{MiFr05} R. Miron and C. Frigioiu, Algebras Groups Geom. \textbf{22}%
, 151 (2005).

\bibitem{MHSS} R. Miron, D. Hrimiuc, H. Shimada and V. S. Sabau, \textit{The
Geometry of Hamilton and Lagrange Spaces}, Kluwer Acad. Publ., Dordrecht;
Boston (2001).
\end{thebibliography}
\end{document}